\documentclass[12pt]{article}
\usepackage{amsmath,amsfonts,amssymb,amsthm}
\def\noi{\noindent}

\def\pf{\noi{\bf Proof.\ \,}}
\def\eop{{$\square$}}

\def\em{\it }        

\def\a{\alpha}
\def\b{\beta}
\def\g{\gamma}
\def\d{\delta}
\def\l{\lambda}
\def\o{\omega}

\def\th{\theta}
\def\vep{\varepsilon}

\def\CC{{\mathbb C}}
\def\FF{{\mathbb F}}

\def\MM{{\mathbb M}}
\def\QQ{{\mathbb Q}}

\def\TT{{\mathbb T}}
\def\ZZ{{\mathbb Z}}

\def\la{\langle}
\def\ra{\rangle}
\def\<{\langle}
\def\>{\rangle}

\def\bs{\it}            

\def\dim{{\bs dim}}

\def\exp{{\bs exp}}
\def\det{{\bs det}}

\def\l{{\lambda}}

\def\wt{{\rm wt}}

\def\half{{1 \over 2}}
\def\fourth{{1 \over 4}}

\def\sixteenth{{1 \over 16}}

\def\vac{\hbox{\bf 1}} 

\def\vnat{V^\natural}

\def\veh{V_{E_8}}

\def\veh{V_{E_8}}

\def\veehp{{V_{EE_8}^+}}

\def\labtt#1{\label {#1}}

\def\refpp#1{(\ref {#1})}

\begin{document}

\newtheorem{thm}{Theorem}[section]
\newtheorem{prop}[thm]{Proposition}
\newtheorem{lem}[thm]{Lemma}
\newtheorem{rem}[thm]{Remark}
\newtheorem{coro}[thm]{Corollary}
\newtheorem{conj}[thm]{Conjecture}
\newtheorem{de}[thm]{Definition}
\newtheorem{hyp}[thm]{Hypothesis}

\newtheorem{nota}[thm]{Notation}
\newtheorem{ex}[thm]{Example}
\newtheorem{proc}[thm]{Procedure}

\begin{center}\end{center}

\centerline{ January 15, 2012 }
\begin{center}
{\Large \bf  Integral forms in vertex operator algebras which are invariant under finite groups }

\vskip 1cm

Chongying Dong

Department of Mathematics,

University of California,

Santa Cruz, CA 95064 USA

{\tt dong@ucsc.edu}

\vskip 0.5cm
and
\vskip 0.5cm

Robert L. Griess Jr.

Department of Mathematics,

University of
Michigan,

Ann Arbor, MI 48109-1043  USA

{\tt rlg@umich.edu}

\end{center}
\begin{abstract}
For certain vertex operator algebras (e.g., lattice type) and given finite group of automorphisms,
we prove existence of a positive definite integral form invariant under the group.
Applications include an integral form in the Moonshine VOA which is invariant under the Monster, and examples in other lattice type VOAs.
\end{abstract}

\tableofcontents

\section{Introduction}
\def\theequation{1.\arabic{equation}}
\setcounter{equation}{0}

Given a finite subgroup $G$ of $Aut(V)$, where $V$ is a vertex operator algebras (VOA), one would like to study possible $G$-invariant integral forms of $V$.  We give methods to do so in a few cases.  An application is a Monster-invariant form within the Moonshine VOA \cite{FLM}.
(The Monster  is the largest sporadic finite simple group \cite{grfg}.)

An integral form of a VOA
could be considered an analogue of the integral structure spanned by a Chevalley basis in a Lie algebra or of an integral structure in an enveloping algebra.
An integral form of a VOA allows one to create VOAs over other commutative rings of coefficients, e.g., finite fields, algebraic number fields.  When such a form is invariant under a finite group, we get an infinite series of integral representations and actions of the group as automorphisms of finite rank algebras over $\ZZ$ (such algebras are typically nonassociative).  The form may be reduced modulo primes to create modular representations of that group.

An integral form $R$ of a vertex operator algebra $V$ with a nondegenerate symmetric invariant bilinear form $(\cdot,\cdot)$ is a vertex algebra over $\ZZ$ satisfying $(R,R)\le  \QQ.$
Our main examples are lattice vertex operator algebras \cite{B1986}, \cite{FLM} and their orbifolds and subalgebras. There is also a positive definite hermitian form on lattice vertex operator algebras \cite{B1986} although the form is not invariant in the sense of \cite{FHL}.  Hermitian forms invariant by a finite group also give integral representations.
Some of our examples satisfy the stronger condition $(R,R)\le \ZZ$.

\section{Preliminaries}
\def\theequation{2.\arabic{equation}}
\setcounter{equation}{0}

\begin{de}\labtt{ratform}   Let $\FF$ be a field of characteristic 0.  A rational form of an $\FF$-vector space is a rational subspace for which any basis is an $\FF$-basis of the complex subspace.
\end{de}

\begin{de}\labtt{intform}  An integral form in a rational vector $U$ space of dimension $m$ is a free abelian subgroup of $(U,+)$ of rank $m$ (so a basis of a lattice is a basis of $U$ as a vector space).
If $U$ has a symmetric bilinear form,
the integral form is an integral (resp., rational) lattice if  the integral form is an integral (resp., rational) lattice with respect to the form.

Let $\FF$ be a field of characteristic 0.
An integral form, etc. for an $\FF$-vector space is an integral form, etc. for
a rational form of the $\FF$-vector space.
\end{de}

For a VOA, we use the usual symbols $\vac$ for vacuum element and $\o$ for the (principal) Virasoro element.

\begin{de}\labtt{intformvoa}  Suppose that $V$ is a VOA (over the complex numbers) with a nondegenerate symmetric bilinear form.
An integral VOA form (abbreviated IVOA) for $V$ is an abelian subgroup $R$ of $(V, +)$ such that $R$ is a VA over $\ZZ,$ there exists a positive integer $s$ so that $s\omega\in R,$ for each $n$, $R_n:=R\cap V_n$ is an integral form of $V_n$, $(R,R)\le \QQ$.
Since $R$ is a VOA,  $\vac\in R$.   For each degree $n$,  $R_n$ has finite rank, whence there is an integer $d(n)>0$ so that $d(n)\cdot (R_n,R_n) \le \ZZ$.
\end{de}

\begin{rem}\labtt{aboutR}  (i) If $\dim (V_0)=1,$  the invariant property of the bilinear form  requires that for any $u,v\in R_m,$
 $\sum_{n\geq 0}\frac{1}{n!}(L(1)^nu)_{2\wt u -2n-1}v\in \frac{1}{d(m)}\ZZ \vac.$

(ii) If $R$ is generated over $\ZZ$ by a set of homogeneous elements, then $R=\bigoplus_{n \in \ZZ} R_n$ is homogeneous.
\end{rem}

\begin{thm} \labtt{invariantintform}
We use the notations of \refpp{intformvoa}.
We assume that the finitely generated
$V$ is generated over $\CC$
by
$\tilde R := R_0+R_1+\cdots +R_d$ and that the finite group
$G$ leaves invariant the rational vector space spanned by $\tilde R$.

There exists an integral form $S$ of $V$ such that

(a) $S\le R$ and for all $n$, $|R_n : S_n|$ is finite (where $S_n:=S\cap V_n$);

(b) $S$ is $G$-invariant;

\end{thm}

\pf    For all $n$, define $U_n$ to be the $\QQ$-span of $R_n$.

Since $G$ is finite,
the $G$-invariant abelian group $S:=\bigcap_{g\in G} g{R}$
has the property that for all $n$,
$|S_n: R_n|$ is finite.  Therefore $S_n$ is a lattice in $U_n$.

Since $G\le Aut(V)$, $gR$ is an integral form of $V$, for all $g \in G$, and so $S$ is also an integral form.
\eop
\section{Integral forms of lattice vertex algebras}
\def\theequation{3.\arabic{equation}}
\setcounter{equation}{0}

Let $L$ be a positive definite even lattice with a basis  $\{\gamma_1,\dots ,\gamma_d\}$ and $L^{*}$ the dual lattice of $L.$ We will follow the setting of \cite{FLM} for lattice vertex operator algebra $V_L$ and the Heisenberg vertex operator subalgebra $M(1)$ of $V_L.$

For any
$\alpha\in L^{*}$ we set
$$E^-(-\alpha,z):=\exp\left(\sum_{n>0}\frac{\alpha(-n)}{n}z^n\right)=\sum_{n\geq 0}s_{\alpha,n}z^n.$$
Let $M(1)_{\ZZ}$ be the ${\ZZ}$-span of  $s_{\alpha_1,n_1}s_{\alpha_2,n_2}\cdots s_{\alpha_k,n_k}$
for $\alpha_i\in L$ and $n_i\geq 0.$ Then $M(1)_{\ZZ}$ is a subset of $M(1)$ which is a vertex operator subalgebra of $V_L.$
We will prove that  $M(1)_{\ZZ}$ is, in fact, an IVOA.

We will denote the set of partitions of postive integers by ${\cal P}.$
For each $\alpha\in L^{*}$ and $\lambda=(\lambda_1,\lambda_2,...,\lambda_k)\in{\cal P}$ with $\lambda_1\geq \lambda_2\geq \cdots \geq \lambda_k>0$
we define vectors by the $k\times k$ determinant
$$s_{\lambda}(\alpha)=\det(s_{\alpha,\lambda_i+j-i})$$
where $s_{\alpha, n}$ is understood to be zero if $n<0.$
Then $s_{\lambda}(\alpha)\in M(1)_{\ZZ}$ for any $\alpha\in L$ and $\lambda\in {\cal P}.$

\begin{lem}\labtt{firstfew}  We have
$E^-(\a, z) = $
\smallskip

$(1z^0 + \a (-1)z^1 +\half \a (-1)^2  z^2 +\frac 16 \a (-1)^3 z^3 + \cdots ) \times $

\smallskip

$(1z^0+ \frac 12 \a (-2)z^2 + \frac 12  \frac 1 {2^2} \a (-2)^2 z^4 + \frac 16 \frac 1{2^3} \a (-2)^3 z^6 + \cdots ) \times $

\smallskip

$(1 z^0 + \frac 13 \a (-3) z^3 + \frac 12 \frac 1{3^2} \a (-3)^2 z^6 + \frac 16 \frac 1{3^3} \a (-3)^3 z^9  + \cdots ) \times  \cdots $.

\bigskip

\noindent
The  series $E^-(\a, z)$ begins:
$1z^0 \ + \ \a(-1) z^1 \ + \
(\half \a (-2) +\half \a (-1)^2 ) z^2 \ + \ $

\smallskip

\noindent
$(\frac 13 \a (-3) + \half \a (-1)\a(-2) + \frac 16 \a (-1)^3) z^3 \ + \  \cdots .$

\smallskip

It follows that $s_{\a , 0}=1, s_{\a ,1}= \a (-1), s_{\a, 2}=\half \a (-2) +\half \a (-1)^2$ and $s_{\a, 3}= \frac 13 \a (-3) + \half \a (-1)\a(-2) + \frac 16 \a (-1)^3$.

\end{lem}

\begin{lem}\labtt{ld1} $M(1)_{\ZZ}$ has a ${\ZZ}$-basis
$$s_{\alpha_1,n_1}s_{\alpha_2,n_2}\cdots s_{\alpha_k,n_k}$$
where $\alpha_i\in \{\gamma_1,\dots \gamma_d\},$ $n_1\geq n_2\cdots \geq n_k$ and $k\geq 0.$
\end{lem}

\pf  Since $\{\gamma_1,\dots ,\gamma_d\}$ is a basis of $L$ we see that $M(1)_{\ZZ}$ is spanned
by
$$s_{\alpha_1,n_1}s_{\alpha_2,n_2}\cdots s_{\alpha_k,n_k}$$
for $\alpha_i\in \{\pm \gamma_1,\dots , \pm\gamma_d\}$ and $k\geq 0.$
 It is enough to show that
 each
 $s_{-\gamma_i,n}$ lies in the ${\ZZ}$-span
 of $$s_{\alpha_1,n_1}s_{\alpha_2,n_2}\cdots s_{\alpha_k,n_k}$$
where $\alpha_i\in \{\gamma_1,\dots \gamma_d\}$ and $k\geq 0.$
Note that $E^-(-\alpha,z)E^-(\alpha,z)=1$ for $\alpha\in L.$
Then for any $n\geq 0$
$$\sum_{i=0}^ns_{\a,i}s_{-\a,n-i}=\delta_{n,0}.$$
Since each $s_{\a , 0}=1$, it
follows that $s_{-\a,n}$ is a $\ZZ$-linear combination of products of $s_{\a,m}$'s.
The ${\ZZ}$-linear independence of these vectors follows from the ${\CC}$-linear independence.
\eop

\begin{thm}\labtt{td2} Let $(V_L)_{\ZZ}$ be the $\ZZ$-span of
$$s_{\alpha_1,n_1}s_{\alpha_2,n_2}\cdots s_{\alpha_k,n_k}e^{\alpha}$$
where $\alpha_i\in \{\gamma_1,\dots \gamma_d\},$ $n_1\geq n_2\cdots \geq n_k,$
$k\geq 0$ and $\alpha\in L.$ Then $(V_L)_{\ZZ}$ is a IVOA generated by $e^{\pm \gamma_i}$ for $i=1,\dots ,d.$  The natural bilinear form restricted to $(V_L)_{\ZZ}$ is integral. Moreover, $M(1)_{\ZZ}$ is a sub IVOA of $(V_L)_{\ZZ}.$
\end{thm}

\pf  From lemma \ref{ld1} and the construction of $V_L,$ the spanning set of $(V_L)_{\ZZ}$ forms a ${\ZZ}$-basis of $V_L.$ We first prove that $(V_L)_{\ZZ}$ is a vertex algebra over $\ZZ.$
That is, we need to show that for basis element $u,v$ and $n\in \ZZ,$ $u_nv\in (V_L)_{\ZZ}.$

Let  $\alpha,\beta\in L.$ Let $\alpha_1,\dots ,\alpha_k,
\beta_1,\dots ,\beta_l\in L$ subject to the conditions $\alpha_1+\cdots+ \alpha_k=\alpha$ and
$\beta_1+\cdots+\beta_l=\beta.$ Set
$$A:=\exp\left(\sum_{i=1}^k\sum_{n\geq1}\frac{\alpha_i(-n)}{n}w_i^n\right)e^{\alpha}\in V_L[[w_1,\cdots,w_k]]$$
$$B:=\exp\left(\sum_{j=1}^l\sum_{n\geq1}\frac{\beta_j(-n)}{n}x_j^n\right)e^{\beta}\in V_L[[x_1,\cdots,x_l]]$$
Then the coefficients in the formal power series $A$ and $B$ span $M(1)\otimes e^{\alpha}$ and
$M(1)\otimes e^{\beta}$ respectively as $k,l,\alpha_i$ and $\beta_j$ vary subject to the conditions
(cf. \cite{FLM}). It is good enough to prove that the coefficients of $Y(A,z)B$ lie in $(V_L)_{\ZZ}$ again.

It follows from \cite{FLM} that
$$Y(A,z)B=\mbox{$\circ\atop\circ$}Y(e^{\alpha_1},z+w_1)\cdots Y(e^{\alpha_k},z+w_k)\ \times $$
$$Y(e^{\beta_1},x_1)\cdots Y(e^{\beta_l},x_l)\mbox{$\circ\atop\circ$}{\bf 1}\prod_{1\leq i\leq k, 1\leq j\leq l}
(z+w_i-x_j)^{(\alpha_i,\beta_j)}.$$
Since the coefficients of $\prod_{1\leq i\leq k, 1\leq j\leq l}
(z+w_i-x_j)^{(\alpha_i,\beta_j)}$ are integers, it is good enough to show that the coefficients
of $\exp\left(\sum_{n>0}\frac{\alpha_i(-n)}{n}(z+w_i)^n\right)$ preserve $(V_L)_{\ZZ}.$  Again,  this is clear
from the definition of $(V_L)_{\ZZ}$ and Lemma \ref{ld1}.

It is well known from the construction of $V_L$ \cite{FLM} that $V_L$ is generated by $e^{\pm \gamma_i}$ for $i=1,\dots .,d.$  Let $V$ be a vertex operator algebra $V$ and $u,v\in V.$ Then for any $n\in {\ZZ}$ the vertex operator
$$Y(u_nv,z)=Res_{z_1}\{(z_1-z)^nY(u,z_1)Y(v,z)-(-z+z_1)Y(v,z)Y(u,z_1)\}$$
involves only  integer linear combinations of products  $u_sv_t$ and
$v_su_t$ for $s,t\in {\ZZ}.$  As a result, $(V_L)_{\ZZ}$ is generated by $e^{\pm \gamma_i}$ for $i=1,\dots , d.$

Since $\dim (V_L)_0=1$ and $L(1)(V_L)_1=0$, there is a unique non-degenerate symmetric invariant bilinear form on $V_L$ \cite{L}, Prop. 3.1.    It is characterized by the conditions below.

(1) $(e^\alpha,e^\beta)=\delta_{\alpha,-\beta}$ for $\alpha,\beta\in L$
where we have arranged  for $e^{\alpha}e^{\beta}=\epsilon(\alpha,\beta)e^{\alpha+\beta}$
for a bimultiplicative form
$\epsilon: L\times L\to \<\pm 1\>$ such that
$\epsilon(\alpha,\alpha)=(-1)^{(\alpha,\alpha)/2}.$

(2) $(\alpha(n)u,v)=-(u,\alpha(-n)v)$ for all  $u,v\in (V_L)_{\ZZ},$    $\alpha\in L$ and $n\in {\ZZ}.$

We show that $((V_L)_{\ZZ}, (V_L)_{\ZZ})\leq \ZZ.$  Let $\alpha_1,\dots ,\alpha_k\in L^{*}$, $\alpha,
\beta_1,\dots ,\beta_l,\beta\in L$ and $A,B$ be as before.
Then we have
\begin{equation}\label{3.1}(A,B)=(e^{\alpha}, \exp\left(-\sum_{i=1}^k\sum_{n\geq1}\frac{\alpha_i(n)}{n}w_i^n\right)\exp\left(\sum_{j=1}^l\sum_{n\geq1}\frac{\beta_j(-n)}{n}x_j^n\right)e^\beta)
\end{equation}
$$=\prod_{i,j}(1-w_ix_j)^{(\alpha_i,\beta_j)}\delta_{\alpha,-\beta}.$$
That is,
$$\sum(s_{\alpha_1,m_1}\cdots s_{\alpha_k,m_k}e^{\alpha},s_{\beta_1,n_1}\cdots s_{\beta_l,n_l}e^{\beta})w_1^{m_1}\cdots w_k^{m_k}x_1^{n_1}\cdots x_l^{n_l}$$
$$=\prod_{i,j}(1-w_ix_j)^{(\alpha_i,\beta_j)}\delta_{\alpha,-\beta}.$$

Since $M(1)$ is a sub VOA of $V_L,$ it follows immediately that $M(1)_{\ZZ}$ is a sub IVOA of $(V_L)_{\ZZ}.$
The proof is complete. \eop

\begin{rem}\labtt{hermitian}
There is also a unique positive definite hermitian form $(\cdot\mid \cdot)$ on $V_L$ \cite{B1986} such that

(1) $(e^\alpha\mid e^\beta)=\delta_{\alpha,\beta}$ for $\alpha,\beta\in L$;

(2) $(\alpha(n)u\mid v)=(u\mid \alpha(-n)v)$ for all  $u,v\in (V_L)_{\ZZ},$ $\alpha\in L$ and $n\in {\ZZ}.$

Note that the positive definite form $( \cdot  \mid  \cdot  )$ on $V_L$ is not invariant in the sense of \cite{FHL}  but it is  useful.

As in the discussion for the bilinear form, we have
$$(A\mid B)=(e^{\alpha}\mid  \exp\left(\sum_{i=1}^k\sum_{n\geq1}\frac{\alpha_i(n)}{n}w_i^n\right)\exp\left(\sum_{j=1}^l\sum_{n\geq1}\frac{\beta_j(-n)}{n}x_j^n\right)e^\beta)
$$
$$=\prod_{i,j}(1-w_ix_j)^{-(\alpha_i,\beta_j)}\delta_{\alpha,\beta}.$$
That is,
$$\sum(s_{\alpha_1,m_1}\cdots s_{\alpha_k,m_k}e^{\alpha}\mid s_{\beta_1,n_1}\cdots s_{\beta_l,n_l}e^{\beta})w_1^{m_1}\cdots w_k^{m_k}x_1^{n_1}\cdots x_l^{n_l}$$
$$=\prod_{i,j}(1-w_ix_j)^{-(\alpha_i,\beta_j)}\delta_{\alpha,\beta}.$$
\end{rem}

We will next determine the dual lattice of the integral lattice $(V_L)_{\ZZ}$ in $V_L$ with respect to the forms $(\ , )$ and $(\, \mid\,)$.

\begin{nota}\labtt{nota1}   Let $\{\beta_1,\dots ,\beta_d\}$ be the  basis of the dual lattice $L^{*}$ which is dual to
the basis
$\{\gamma_1,... , \gamma_d\}$  of $L.$
 Let $U$ be the ${\ZZ}$-span of
the ${\ZZ}$-linearly independent vectors
$$s_{\lambda_1}(\beta_1)\cdots s_{\lambda_d}(\beta_d)\otimes e^{\a}$$
for $\lambda_1,...,\lambda_d\in {\cal P}$ and $\a\in L.$
If $L$ is self dual, then $(V_L)_{\ZZ}=U.$
\end{nota}

\begin{prop}\labtt{uisdual}
(i) The 
set
$$\{s_{\lambda_1}(\gamma_1)\cdots s_{\lambda_d}(\gamma_d)\otimes e^{\a}\mid \lambda_1,...,\lambda_d\in {\cal P}, \a\in L\}$$
is a $\ZZ$-basis of $V_L)_{\ZZ}$.

(ii) $U$ is  the dual of $(V_L)_{\ZZ}$ in $V_L$ with respect to both the bilinear form $( , )$ and the hermitian form $(\ \mid \ ).$

(iii) The set
$$\{s_{\lambda_1}(-\beta_1)\cdots s_{\lambda_d}(-\beta_d)\otimes e^{-\a}\mid \lambda_1,...,\lambda_d\in {\cal P},\a\in L\}$$
is the dual basis
with respect to $( ,)$ and
the set
$$\{s_{\lambda_1}(\beta_1)\cdots s_{\lambda_d}(\beta_d)\otimes e^{\a}\mid \lambda_1,...,\lambda_d\in {\cal P},\a\in L\}$$
is the dual basis
with respect to $(\ \mid \ ).$
In particular, if $L$ is self dual, then $(V_L)_{\ZZ}$ is also self dual with respect to both the bilinear form $(\, ,\, )$ and the hermitian form $(\,\mid\,  ).$
\end{prop}
\pf  (i, ii) We deal only with the hermitian form and the proof for the bilinear form is similar.

We first prove that
$$(s_{\lambda_1}(\beta_1)\cdots s_{\lambda_d}(\beta_d)\otimes e^{-\a} \mid  s_{\lambda_1'}(\gamma_1)\cdots s_{\lambda_d'}(\gamma_d)\otimes e^{\beta})=\delta_{(\lambda_1,...,\lambda_d),(\lambda_1',...,\lambda_d')}\delta_{\alpha,\beta}$$
for $\lambda_i,\lambda_i'\in {\cal P}$ and $\alpha,\beta\in L.$ Since
 $(\beta_i,\gamma_j)=\delta_{i,j}.$ It is enough to show that
 $$(s_{\lambda}(\beta_i) \mid s_{\lambda'}(\gamma_i))=\delta_{\lambda,\lambda'}$$
 for $\lambda,\lambda'\in {\cal P}$ and $i=1,...,d.$ But this is clear from the bilinear form formula
 $(\ref{3.1})$ (see Formula (4.8) of \cite{M}).

 We next show that the set of
 $$s_{\lambda_1}(\gamma_1)\cdots s_{\lambda_d}(\gamma_d)$$
for $\lambda_1,...,\lambda_d\in {\cal P}$ form a basis of $M(1)_{\ZZ}.$
Note that for $\alpha_1,...,\alpha_k\in L$ and $n_1,...,n_l\geq 0,$
$$(s_{\lambda_1}(\beta_1)\cdots s_{\lambda_d}(\beta_d) \mid s_{\a_1,n_1}\cdots s_{\a_k,n_k})\in \ZZ$$
for $\lambda_i\in {\cal P}.$ This implies that  $s_{\a_1,n_1}\cdots s_{\a_k,n_k}$ is a $\ZZ$-linear combination of $s_{\lambda_1}(\gamma_1)\cdots s_{\lambda_d}(\gamma_d)$ for $\lambda_i\in {\cal P}.$ Consequently, $$s_{\lambda_1}(\gamma_1)\cdots s_{\lambda_d}(\gamma_d)$$
for $\lambda_1,...,\lambda_d\in {\cal P}$ form a basis of $M(1)_{\ZZ}$ and
$$s_{\lambda_1}(\gamma_1)\cdots s_{\lambda_d}(\gamma_d)\otimes e^{\a}$$
for $\lambda_1,...,\lambda_d\in {\cal P}$ and $\a\in L$ form a basis of $(V_L)_{\ZZ}.$
So, (i) and (ii) are proved.
The rest of the proposition is clear.
\eop

\begin{rem}
If $L$ is self dual, then $(V_L)_{\ZZ}$ is
a positive definite self dual lattice by Proposition \ref{uisdual}. This was pointed out  in \cite{B1986} already.
\end{rem}

\begin{rem}\labtt{discgroup2} Let $R=(V_L)_{\ZZ}.$ It seems that finding the group structure of $U_n/R_n$ \refpp{nota1} for an arbitrary even lattice $L$ will be very complicated. One can easily see that $U_1/R_1$ is isomorphic to $L^{*}/L.$
\end{rem}

\subsection{About the automorphism group of $V_L$ and a torus normalizer} \labtt{liftO(L)}

\begin{nota}\labtt{elemabel}
In this article, we use the standard notation $p^n$ for an elementary abelian $p$-group of order $p^n$, where $p$ is a prime.   Context will indicate when $p^n$ means a group.
\end{nota}

We now  review the automorphism group of $V_L$ and certain other groups which act on $(V_L)_{\ZZ}$ preserving either the bilinear form or the hermitian form.
For this purpose we need to fix a bimultiplicative function
$\epsilon:L\times L\to \langle \pm 1\rangle$ such that $c(\alpha,\beta):=\epsilon(\a,\b)\epsilon(\b,\a)$ equals $(-1)^{(\a,\b)}$
for $\a, \b\in L.$ Let $\hat L$ be the central extension of $L$ by $\langle \pm 1\rangle$ with commutator map $c$ and we denote the map from $\hat L$ to $L$ by $\bar{}$. Then any automorphism $\sigma$ of $\hat L$
induces an automorphism $\bar \sigma$ of $L.$   Let $O(L)$ be the  isometry group of $L$ and
$O(\hat L):=\{\sigma\in Aut(\hat L) \mid \bar \sigma\in O(L)\}.$
The above group embeds in $Aut(V_L)$ and permutes the set $\{ \pm e^{\a} \mid \a \in L   \}$.   This subgroup of $Aut(V_L)$ is often denoted $\widetilde{O(L)}$ in the literature.
It has the shape $\widetilde{O(L)} \cong \ZZ^{rank(L)}.O(L)$.

Let $\TT$
be the torus obtained by exponentiating $a_0$ where $a$ ranges over the standard Cartan subalgebra of $(V_L)_1$.
Then $\TT \cong H/L^{*}$ where $H=\CC\otimes_{\ZZ}L$ and
 the normalizer $N_L$ of $\TT$ in $Aut(V_L)$  is a group of the form
$N_L \cong \TT.O(\hat L).$ 
We let $N$ be the normal subgroup of $Aut(V_L)$ generated by $e^{v_0}$ for
$v\in (V_L)_1.$ Then $Aut(V_L)=N\, N_L$, product of subgroups,  and $N$ is the connected component of the identity in the algebraic group  $Aut(V_L)$.  See \cite{DN} and \cite{dgag} for background.

There  exists a standard lift, called $\theta$,  of the $-1$ isometry of $L$ to an automorphism of order 2 of $V_L$.  The standard lift interchanges $e^{\a}$ and $e^{-\a}$.    See \cite{FLM} and an appendix of
\cite{ghframes} for a general discussion of lifts.   We have $C_{N_L}(\theta ) \cong 2^{rank(L)}.O(L).$

\begin{lem}\labtt{invariantivoa}  Suppose that $L$ is an even lattice.
We use the IVOA $(V_L)_{\ZZ}$ of \refpp{td2}.
If $S$ is any subgroup of $\widetilde {O(L)}$,  the fixed point subVOA $V^S$
has the IVOA $(V_L)_{\ZZ}^S$.
Note that $N_{\widetilde {O(L)}}(S)/S$
acts as automorphisms of the IVOA $(V_L)_{\ZZ}^S$.
If $M$ is a sublattice so that $S$ fixes the subVOA
$V_M$, then $V_M \cap (V_L)_{\ZZ}^S$ is an IVOA in $V_M^S$.
\end{lem}

\begin{nota}\labtt{ivoa+}
A special case of $S$ in \refpp{invariantivoa}   is  the fixed points of $\theta $ \refpp{liftO(L)} in $(V_L)_{\ZZ}$, denoted  by   $(V_L)_{\ZZ}^+.$  Then
$(V_L)_{\ZZ}^+$ is an IVOA in $V_L^+.$
\end{nota}

\begin{rem}\labtt{twoforms}  We discuss the hermitian form more.  We  define the anti-linear  involution $\sigma$ from $V_L$ to $V_L$ such
 that $\sigma (au)=\bar{a}\theta(u)$ for $a\in \CC$ and $u\in (V_L)_{\ZZ}$ where $\bar{a}$ is the complex conjugate of $a$. It is easy to see that $(u \mid v)=(u,\sigma(v))$ for $u,v\in V_L.$   Note that $N_L$ does not preserve the hermitian form
but the subgroup $\widetilde {O(L)}$ does since it preserves the set
$\{\pm e^{\a} \mid \a \in L \}$, which  generates the $\QQ$-VOA $\QQ \otimes (V_L)_{\ZZ}$.
 \end{rem}

\begin{lem} \labtt{tensor}  If $A, B$ are IVOA forms on the respective VOAs $V, W$, then $A\otimes B$ is an IVOA form on $V\otimes W$ (graded tensor product of VOAs; see \cite{FHL}).
\end{lem}

\subsection{Initial pieces of the standard integral form for $\veh$}

Let $L$ be any even integral lattice. Recall the bilinear form $( \, ,)$ and the hermitian form $(\cdot \mid \cdot )$ on
$V_L$ \refpp{hermitian}.
We use the usual notation for
 $H=\CC\otimes_{\ZZ} L$ and the subspaces $H_{-n}=H\otimes t^{-n}$ of  $V_L$.

The next lemma recalls a few calculations for our use.

\begin{lem}\labtt{someip}
(i)
For $p, q, r, s \in H_{-1}$,

$(pq\mid rs)=(p\mid r)(q \mid s)+ (p \mid s)(q \mid r)$.

(ii)
The elements listed in \refpp{firstfew} at the $z^2$-coefficient have the following inner products (the right side means inner product in $L$).

$(\a (-2) \mid  \b (-1)^2)=0$;

$(\a (-1)\b (-1) \mid  \g (-1)\d (-1))= (\a , \g)(\b, \d) + (\a, \d)(\b, \g)$;

$(\a (-2)\mid \b (-2))=2 (\a, \b )$.

$(\a (-1)^2\mid\a (-1)^2)=2 (\a, \a )^2$, which is 8 if $\a$ is a root.

$(\a (-2)\mid\a (-2))=2 (\a, \a )$, which is 4 if $\a$ is a root.

So if $\a$ is a root,  $\half \a (-2) +\half \a (-1)^2 $ has norm $1+2=3$.

$( \half \a (-2) +\half \a (-1)^2\mid  \half \b (-2) +\half \b (-1)^2)= \half (\a, \b ) + \half (\a, \b)^2  $.
\end{lem}

\begin{nota}  \labtt{E8case}   For the remainder of this section, $L:=E_8$.  We shall study the positive definite lattice $R_2$ with the hermitian form $( \cdot \mid \cdot  )$.
\end{nota}

The lattice $R_2$ generated by these elements is integral and odd (since norm 3 occurs).  Observe that
$R_2$ contains $\half (\a + \b) (-2) +\half (\a + \b)(-1)^2  - (\half \a (-2) +\half \a (-1)^2) -
(\half \b (-2) +\half \b (-1)^2 ) = \a (-1)\b (-1)$, which has odd norm if $(\a, \b)$ is odd.  Note that
$\a (-1)\b (-1)$ is in $ span_{\ZZ} \{ \half \a (-1)^2 \mid \a \in L, (\a, \a)=2 \}$.

\begin{prop}\labtt{latticedegree1and2}
We use the notation \refpp{intformvoa} for the IVOA created in \refpp{td2} and we use  \refpp{firstfew}.

{\bf Degree 1. }
$R_1$ is isometric to $E_8 \perp Q$, where $Q$ is the square lattice of rank 240 spanned by an orthonormal basis (i.e. $det(Q)=1$).  Therefore, $R_1$ is unimodular of rank 248, but not even.

{\bf Degree 2. }  $R_2$ is isometric to $J \perp (E_8 \otimes Q) \perp S$, where $Q$ is the unimodular square lattice of rank 240, $S$ is a unimodular square lattice of rank 2160 and $J$ is a rank 44 determinant 1 odd lattice.
\end{prop}
\pf
The statement about degree 1 is easy.

The lattice  $R_2$ has a basis
$$\{e^{\alpha}\mid \alpha\in E_8, (\a |\a)=4 \}\cup \{\gamma_i(-1)\otimes e^{\alpha}\mid 1\leq i\leq 8, \alpha\in E_8, (\a | \a)=2 \}$$
$$\cup \{s_{\gamma_i,2}, \ s_{\gamma_i,1}s_{\gamma_j,1}\mid 1\leq i,j\leq 8, i\leq j\}.$$
Here, the $\g_i$ form  a basis of the lattice $E_8$ as before.
Note that $s_{\alpha,2}=\frac{1}{2}(\alpha(-1)^2+\alpha(-2))$ and $s_{\alpha,1}=\alpha(-1)$.  The three sets displayed are pairwise orthogonal.

Now to describe $R_2$.  It is isometric to $S \perp (E_8 \otimes Q) \perp J$,
corresponding to the sublattices spanned by the respective
three parts of the basis described  in the last paragraph. Here,  $Q$ is as described in Degree 1, so is square, whence $E_8 \otimes Q$ is the orthogonal direct sum of 240 copies of $E_8$.

The description of the rank 44 lattice $J$ is more complicated.
Write $J_1$ for $J\cap S^2(H_{-1})$ and write $J_2$ for $J \cap H_{-2}$.
From \refpp{sym2lattice} and $det(E_8)=1$, we conclude that the discriminant group of  $J_1$ is elementary abelian of order $2^8$.   It is clear that $J_2$ is spanned by all $\a (-2)$ and so $J_2 \cong \sqrt 2 E_8$ and $J_2$ has discriminant group which is elementary abelian of order $2^8$.

The lattice $J$ represents a gluing of $J_1$ and $J_2$ (with glue vectors $s_{\a, 2} = \half \a (-1)^2 +\half \a (-2)$).  So $J$ contains $J_1\perp J_2$ with index $2^8$, whence $det(J)=1$.   The minimum norm of $J$ appears to be 3.
\eop

\def\cvcch{cvcc\frac 12}
\section{About a given rational form and $\cvcch$}

\begin{nota}\labtt{cvcch} We abbreviate conformal vector of central charge $\half$ by $\cvcch$.  The subVOA generated by a set of elements $w, \dots$ is denoted $subVOA(w, \dots )$.
\end{nota}

\begin{lem}\labtt{cvcchrat1}
Suppose that $Q$ is a rational form of the VOA $V$ and $w\in Q$ is a $\cvcch$.
Then
the Miyamoto involution $t(w)$ for  $w$ stablilizes $Q$.
\end{lem}
\pf
Since $t(w)$ centralizes $subVOA(w)$ and leaves invariant its highest weight spaces,
$V'(c)$ in $V$, for highest weights $c\in \{0, \frac 12, \frac 1{16} \}$, it leaves invariant the
corresponding respective homogeneous components, $V(c)$.

By hypothesis, $Q$ is a rational form so that $Q$ is invariant under each $w_i$, where $Y(w,z)=\sum_{i\in \ZZ} w_i z^{-i+1}$.

Since $V'(c) =(\cap_{i \ge 2} Ker_V(w_i)) \cap Ker_V(w_1 - c)$ and $0, c $ are rational,
$Q'(c):=V'(c)\cap Q$ is a rational form of $V'(c)$.
Since the $w_i$ leave $Q$ invariant,
$Q(c):=Q\cap V(c)$ is a rational form of $V(c)$, the homogeneous component.

We have $Q=Q(0)\oplus Q(\frac 12 )\oplus Q(\frac 1{16})$ The action of $t(w)$ on
each $V(c) \ge Q(c)$ is a scalar, $\pm 1$ (depending on $c$ and whether $V(\frac 1{16})$ is nonzero).
Therefore, $t(w)$ leaves $Q$ invariant.
\eop

\subsection{An application to $V_{EE_8}^+$}

Up to now, we have seen integral forms which are invariant only under subgroups of $\widetilde {O(L)}$ \refpp{liftO(L)}.
We next prove that there exist IVOA forms in $\veehp$ which are invariant under $Aut(\veehp )\cong O^+(10,2)$.

\begin{nota}\labtt{cvcchee8+}
The set $K$ of $\cvcch$ for $V_{EE_8}^+$ has been listed \cite{gr156}.
Their $\QQ$-span is a $\QQ$-form of $(V_{EE_8}^+)_2$.
The $\ZZ$-span is a rational, non-integral lattice in $(V_{EE_8}^+)_2$, invariant under $Aut(\veehp )$.
Since the degree 2 term generates  the complex VOA $V_{EE_8}^+$  \refpp{vlplusgeneration}, the rational VOA it generates is a rational form of $V_{EE_8}^+$.
\end{nota}

We recall two explicit expressions for $\cvcch$ which often occur in lattice type VOAs.

\begin{nota}\labtt{aa1ee8} A $\cvcch$ of the first kind, or of $AA_1$ type, has the form:

\begin{equation}
e=\sixteenth (t^{-1}\otimes \a)^2 \pm \fourth (e^\a +e^{-\a})
\end{equation}
where $\a$ is a norm 4 vector of $L$ \cite{dmz}.

\medskip

A $\cvcch$ of the second kind, or of $EE_8$ type, has the form:

\begin{equation}
e=\frac 18 q +\frac 1{32}\sum_{\a \in E/\{\pm 1\} } \varphi (\a) (e^\a + e^{-\a})
\end{equation}
where $E$ is a sublattice of $L$, $E\cong EE_8$
(here,  $q$ is a sum $\sum_i  (t^{-1}\otimes u_i)^2$, where $u_i$ is an orthonormal basis of  $\CC\otimes E$;
$\varphi$ is a homomorphism $E\rightarrow \{\pm 1\}$) \cite{dlmn,gr156}.
\end{nota}

\begin{thm}\labtt{integralformvee8+}  $\veehp$  contains an IVOA $R$ which is invariant under
$Aut(\veehp)\cong O^+(10,2)$ and which satisfies $(R,R)\le \ZZ$.
\end{thm}
\pf  
Let  $S$ be an elementary abelian group of order $2^5$ in $Aut(V_L)$ so that $S$ is $2B$-pure \cite{grqs} and take a containment $M\le L$, $M \cong EE_8$, $L \cong E_8$, $V_M$ is $S$-invariant and on $V_L$, the fixed point space for $S$ is  a natural subVOA $V_M^+$ of $V_M$.   By \refpp{invariantivoa},  we get an IVOA $J$  for $\veehp$.
It is integral as a lattice, $(J, J)\le \ZZ$.

The members of $K$ \refpp{cvcchee8+}  are in $W$, the $\QQ$-span of $R$.
The finite group $G \cong O^+(10,2)$ generated by the Miyamoto involutions associated to $K$ leaves both $K$ and $W$ invariant \refpp{cvcchrat1}.  It follows that
$R:=\cap_{g\in G} \  gJ$ is a $G$-invariant IVOA in $\veehp$.  \eop

\begin{rem}\labtt{156lattice}
It would be interesting to know the relationship between $R_2$ (the degree 2 part of $R$, which appears in  the  proof of \refpp{integralformvee8+}) and the rational lattice $span_{\ZZ}(K)$.   For example, how close is one to being a rational multiple of the other?
\end{rem}

\section{Extending an integral form of a subVOA to the full VOA}

\begin{nota}  For subsets $X, Y$ of the VOA $V$, define $X\cdot Y$ to be the $\ZZ$-span of all $x_ky$, $x\in X, y\in Y, k \in \ZZ$ (as usual, $x_ky$ represents the $k$-th composition of the VOA $V$).
\end{nota}

\begin{nota}\labtt{tel}
Suppose that the elementary abelian 2-group $E$ acts on the lattice $L$.
Define $\hat E := Hom(E, \{\pm 1\} )$, the character group of $E$.
If $\l \in \hat E$, let $L_{\l} := \{a\in L\mid ea=\l (e)a \text{ for all }e \in E\}$ be the associated eigenlattice.
The total eigenlattice is $Tel(E,L):=\sum_{\l \in \hat E} L_{\l} $.  The quotient $L/Tel(E,L)$ is annihilated by $|E|$.
\end{nota}

\begin{nota}\labtt{forextend}
We use the following notations and hypotheses for this section.

$V$ is a VOA over the rationals.

$N\le Aut(V)$ has a normal subgroup $E\cong 2^r$ for $r\ge 3$;
 $N$ acts doubly transitively on the nontrivial linear characters of $E$;
 define   $N_{\l} :=Stab_N(\l )$.

 $D$ is a  subgroup of $E$ and $rank(E/D)\ge 2$;
 there exists an IVOA  $A$ in the fixed point subVOA $W:=V^D$
  and $A$ is invariant by $N_N(D)$

We let $X$ be the group of linear characters of $E$ whose kernel contains $D$;
$X$ is an abelian group and $rank(X) = rank(E/D)\ge 2$.

For $\l \in X$, $A_{\l}$ are eigenlattices  \refpp{tel} for $E$, $\l \in X$ (see \refpp{tel} for eigenlattice notation); assume that  $A=\sum_{\l \in X} A_{\l}$.

Set $A'_{1} :=\cap_{g\in N} \, gA_{1 }$ and $A''_{1} :=\sum_{g\in N} \, gA_{1 }$ ($1$ denotes the trivial character).

Fix a linear character $\nu \in X, \nu \ne 1$, then define
$A'_{\nu} :=\cap_{g\in N_{\nu} } \, gA_{\nu }$ and
$A''_{\nu} :=\sum_{g\in N_{\nu} } \, gA_{\nu }$.

For an arbitrary nontrivial  linear character $\l \in \hat E$, define
$A'_{\l} := gA'_{\nu}$ and $A''_{\l} := gA''_{\nu}$,
where $g\in N$ takes $\nu$ to $\l$ (this is independent of choices of $\nu$ and $g$).
Both $A'_{\l} $ and $A''_{\l} $ are $N_{\l}$-invariant and for all $g\in N$,
$gA'_{\l} =A'_{g\l} $ and $gA''_{\l} =A''_{g\l} $.

Set $A':=\sum_{\l \in \hat E} A_{\l}'$
and
$A'':=\sum_{\l \in \hat E} A_{\l}''$ (both $A'$ and $A''$ are $N$-invariant).
\end{nota}

\begin{rem}\labtt{forextend2}
(i)
Note that $A_{\l}$ is not defined unless $\l \in X$, but $A_{\l}'$ and $A_{\l}''$ are defined for all $\l \in \hat E$.

(ii) If $\l \in \hat E$, $g, h \in N$ and $g\nu = \l =  h \nu$, then
$gA_{\nu}$ and $hA_{\nu}$ are integer forms of
the same rational subspace but they are not necessarily equal; $gA_{\nu} \cap hA_{\nu}$ has finite index in each of $gA_{\nu}$ and $hA_{\nu}$.

(iii) The quotient
$A''/A'$ is a torsion abelian group (so that $A'$ and $A''$ have the same rational span).
This follows from (i) and the fact that each $N_{\l}$ is a finite group.
\end{rem}

\begin{lem}\labtt{extend1a}
(i)
$A' \le A'\cdot A' \le  A' \cdot A'' \le A''$.

(ii) $A' \le subVOA(A')\cdot A'' \le A''$.
\end{lem}
\pf
(i)  Since $\vac \in A'$, $A' \le A'\cdot A'$.
We next prove that $A' \cdot A'' \le A''$.

Suppose that the characters $\a, \b$ lie in $X$.   Then
$A_{\a}'\cdot A_{\b} \le  A_{\a}\cdot A_{\b} \le A_{\a \b}\le A_{\a \b}''$.

For each $g\in N_{\b}$, we have
$A_{g\a}' \cdot gA_{\b} \le A_{g(\a \b)}''$.
Fixing $\b$ and summing these containments over all $\a \in X$ and $g \in N_{\b}$, we deduce from double transitivity of $N$ on nonidentity characters and $rank(X)\ge 2$ that
$A' \cdot A_{\b}'' \le A''$.

By using the last containment and transitivity of $N$ on the nontrivial characters in $\hat E$, we conclude that
$A' \cdot A'' \le A''$.

\smallskip

(ii) Let $u\in A', v\in A''.$ By (i) we know that $u_jv\in A'\cdot A''\leq A''$ for all $j\in\ZZ.$ Now assume that $u,v\in V$ such that for any $w\in A''$ and $j\in \ZZ,$ $u_jw \in A''$ and $v_jw\in A''.$ Since $(u_iv)_jw$ for any $i,j\in \ZZ$ is a linear combination of $u_sv_tw, v_pu_qw$ for $s,t,p,q\in\ZZ$ we see that $(u_iv)_jw\in A''.$ As a result, $subVOA(A')\cdot A'' \le A'',$ whence (ii).

\begin{rem}\labtt{extend1a2}
In the notation of \refpp{extend1a}, $subVOA(A')$ is an IVOA and at each degree, the homogeneous summand is a rational lattice of finite rank.   It is unclear if these are integral lattices.   However, we can observe that
the degree 0 summand is just $\ZZ \cdot \vac$.  
We have $\vac \in subVOA(A')$.  There is an integer $m>0$ so that $subVOA(A')\cap \QQ \vac = \ZZ \frac 1m\vac$.  Then, since $\vac _{-1}\vac = \vac$,  $subVOA(A')$ would contain $\frac 1{m^t}\vac$ for all $t\ge 0$.  If $m>1$, this contradicts the hypothesis  that
$subVOA(A') \cap \QQ \vac$ is a rank 1 lattice.
\end{rem}

\begin{thm}\labtt{extend1c}  We use the notation of \refpp{forextend}.   Suppose that $V$ satisfies $dim(V_0)=1$.
There exists an $N$-invariant IVOA whose rational span contains the rational span of $A$.
\end{thm}
\pf
We take $subVOA(A')$, clearly $N$-invariant.
\eop

\subsection{Applications}

We deduce existence of an invariant integral form in a few VOAs of interest.

\begin{thm} \labtt{leechpivoa}  Let $\Lambda$ be the Leech lattice.
$V_{\Lambda}^+$ has an IVOA form which is invariant under $Aut(V_{\Lambda}^+)$ and which furthermore contains $\cvcch$ of both $AA_1$-type and $EE_8$-type.
\end{thm}
\pf
Note that $Aut(V_{\Lambda}^+)$ is just  the centralizer in  $Aut(V_{\Lambda})$ of $\theta$  \refpp{liftO(L)} \cite{ghframes}.   To prove the theorem, just take the fixed points of $\th$ on the IVOA constructed in \refpp{ivoa+}.
\eop

\begin{thm}\labtt{moonivoa}
There exists an IVOA for the Moonshine VOA $\vnat$ which is invariant under $Aut(\vnat )\cong \MM$ and whose rational span contains all $\cvcch$ in $\vnat$.
\end{thm}
\pf
First, we argue that $G:=Aut(\vnat)\cong \MM$ has an invariant rational form in $\vnat$.  This is so, because $Aut(\vnat)$ leaves invariant a rational form in the 196884-dimensional algebra $\vnat_2$ \cite{grfg}, which generates $\vnat$ \cite{FLM}.  Furthermore, this particular rational form contains all $\cvcch$ in $\vnat_2$ since it plainly contains the ones of $AA_1$-type (we identify $V_{\Lambda}^+$ with the fixed points $V^z$ of a $2B$-involution $z$ on $V$)  and the set of all $\cvcch$ is  in bijection with the $2A$-involutions of $\MM$ \cite{miy,hoehn}.

Now to prove existence of an IVOA in this rational form of $\vnat$.
Let $z$ be a $2B$-involution and let $C:=C_G(z)$.

Take any $r\in \{3,4,5\}$.
In $G$, we take a subgroup $E$ which is elementary abelian,
 and $E$
is $2B$-pure of rank $r$ and furthermore satisfies $N_G(E)/E \cong GL(r,2)$.
(With elementary arguments, one can prove existence of such $E$ for $r\le 4$ by looking in $C$ \refpp{2Bpure}.)

Take a subgroup $E$ as in the last paragraph so that $z\in E$ (this is possible since $z$ is in the class $2B$).
We take the IVOA of $V^z$ in \refpp{leechpivoa}, then
we define a $N_G(E)$-invariant IVOA, called $J$,
on $V$ by using \refpp{extend1c}
(in that notation, take $N=N_G(E)$ and $D=\la z \ra$).
It follows that
$\cap_{g\in G}\  gJ$ is a $G$-invariant IVOA and that the
$G$-invariant rational form which it spans contains all $\cvcch$ in $V$.
\eop

\section{Final Remarks}

\subsection{About the dual}

 \begin{lem}\labtt{la1}  Let $R$ be an IVOA.   If $Y(e^{zL(1)}u,z^{-1})v\in R[[z,z^{-1}]]$ for $u,v\in R$ then $(R,R)\le\ZZ$ and $R^*$ is a an $R$-submodule.
 \end{lem}

 \pf  Let $u,v\in R.$ Then
 $$(u,v)=Res_zz^{-1}(Y(u,z), v)=Res_zz^{-1}(\vac, Y(e^{zL(1)}(-z^{-2})^{L(0)}u,z^{-1})v)\in \ZZ.$$
 That is, $(R,R)\le\ZZ.$  Now let $u\in R$ and $w\in R^*.$ We have to show that $Y(u,z)w\in R^*[[z,z^{-1}]].$ This is equivalent
 to showing that for any $v\in R,$ $(Y(u,z)w,v)\in \ZZ[z,z^{-1}]].$ Using the invariant property of the bilinear form we see that
 $$(Y(u,z)w,v)=(w, Y(e^{zL(1)}(-z^{-2})^{L(0)}u,z^{-1})v)\in \ZZ[z,z^{-1}]].$$
 \eop

 \begin{lem}\labtt{la2} If $V=V_L$ then $R^*$ is an $R$-submodule.
 \end{lem}

 \pf By Theorem \ref{td2} we know that $R$ is generated by $e^{\pm\gamma_i}$ for $i=1,...,d.$ From
 the proof of \refpp{la1} it is enough to show that $Y(e^{zL(1)}u,z^{-1})v\in R[[z,z^{-1}]]$  for
 $u=e^{\alpha}$ and $v\in R.$ But this is clear as  $Y(e^{zL(1)}e^{\alpha},z^{-1})v=Y(e^{\alpha},z^{-1})v.$
 \eop

\begin{rem} (i) A vector $v$ in a vertex operator algebra $V$ is called quasi primary vector if $L(1)v=0.$ From the proof of Lemma \ref{la2} we see that if an IVOA $R$ of $V$ 
satisfies $(R,R)\le \ZZ$ and $R$ is generated by the quasi-primary vectors then $R^*$ is an $R$-module.

(ii) If $R^*$ is an $R$-module, $R\cap R^*$ is an IVOA and $(R\cap R^*, R\cap R^*)\le \ZZ$
\end{rem}

\subsection{Examples generated by degree 1}
There are many finite subgroups of Lie groups which have rational-valued characters on the adjoint modules
(see the classification and lists in \cite{grqs}).  Some of these representations are afforded over the rational numbers.  Let $G$ be such a finite subgroup and let $L$ be the root lattice of the Lie algebra.
The  associated representation of $G$ on $Aut(V_L)$, for $L$ the  root lattice, could lead to $G$-invariant IVOAs since the degree 1 component of $V_L$ generates $V_L$ as a VOA.  If so, this procedure would provide an infinite family of integral representations of $G$.

\subsection{Adjusting a rational form to an integral form}

Suppose that we have an IVOA, $R$.
At each degree, $m$, there is an integer $d(m)>0$ so that
$d(m) \cdot (R_m, R_m) \le \ZZ$ \refpp{intformvoa}.

\begin{rem}\labtt{modifyform1}
We can replace the form by the form whose restriction to $R_m$ is  $d(m)$ times the restriction of $( \ , \ )$.   Then we have an IVOA which is an integral lattice, though the modified form may lack properties of the original form, e.g. invariance in the VOA sense.
The modified form is invariant by the subgroup of $Aut(V)$ which leaves $R$ invariant.
\end{rem}
\begin{rem}\labtt{modifyform2}
Another form of interest would be defined on $R_m$ as follows.
On the VOA $V=\bigoplus_{n\in \ZZ}V_n$, the $k$-th product takes $V_i \times V_j$ into $V_{i+j-k-1}$.  The $(m{-}1)$-th product makes $V_m$
a finite dimensional algebra with symmetric bilinear form
$f_m(a,b):=trace(ad(a)\cdot ad(b))$, where $ad(c)$ means the endomorphism
$x\mapsto c_{m-1}x$ of $V_m$.  Since $R$ is an IVOA, the restriction of $f_m$ to $R_m$  is  integer valued.

The form is nondegenerate in some interesting cases, e.g., $m=2$ for the Moonshine VOA.  It is reasonable to ask if it is nondegenerate in general and whether it has interesting applications.   It is invariant under the action of the subgroup of $Aut(V)$ which leaves $R$ invariant.
\end{rem}

\appendix

\section{Appendix: Background}

\begin{lem}\labtt{vlplusgeneration}
Let $L$ be an even lattice which is rootless and is spanned by its norm 4 vectors.  Then $V_L^+$ is generated as a VOA by the elements $e^\a+e^{-\a}$ (in usual notation), where $\a$ ranges over the norm 4 vectors of $L$.
\end{lem}
\pf The argument  in \cite{FLM} for the Leech lattice works here.
\eop

\section{About involutions}

For notation, see \refpp{elemabel}.

\begin{lem}
\labtt{involonfreeabelian}
Suppose that $A$ is a free abelian group and that the involution $t$ acts on $A$.
Define $A^{\vep}$ to be the points of $A$ which are fixed, negated as $\vep = +, -$, respectively.
Suppose that on $A/2A$, the Jordan canonical form for $t$ has exactly $r$ Jordan blocks of size $2\times 2$.   Then $A/(A^+\oplus A^-)\cong 2^r$.
\end{lem}
\pf
Let $B$ be the subgroup defined by $2A\le B\le A$ and $B/2A$ is the fixed point space for the action of $t$ on $A/2A$.  Then $A/B \cong 2^r$.
We shall prove that $B=A^-+A^+$, which implies the result.  

We note that $2A\le (t-1)A + (t+1)A \le A^- + A^+$, whence $B\le A \le \half (A^- + A^+)$.  Take $b \in B$ and write $b=\half (c+d)$ for $c\in A^- , d\in A^+$.   Then $tb=\half (-c+d)$ and so $tb-b=c$, which must be in $2A$ since $b \in B$.   Since $A^-$ is a direct summand of $A$, this means $c\in 2A^-$.  Consideration of $tb+b\in 2A$ proves
$d\in 2A^+$.   We conclude that $B=A^-+A^+$.
\eop

\begin{rem}\labtt{involonfreeabelian2}  The above is well-known.   For more about actions of involutions on free abelian groups, see \cite{gal,bwy,ibw}.
\end{rem}

\section{Appendix: about tensors}

\begin{lem}\labtt{sym2lattice}
Suppose that $K$ is a free abelian group of finite rank $n$ and $t$ is the involution which interchanges tensor factors in $S:=K\otimes K$.   Let $S^+, S^-$ be the points of $S$ which are, respectively,  fixed, negated by $t$.  Then $S/(S^+ + S^-)$ is elementary abelian 2-group of rank $n$.
\end{lem}
\pf
Observe that if $x_1, \dots, x_n$ is a basis of $K$, then the set of all $x_i \otimes x_j$  form a basis of $S$.  This means that the Jordan canonical form of $t$ on $S/2S$ has $n\choose 2$ Jordan blocks of size 2 and $n$ of size 1.  Now use \refpp{involonfreeabelian}.
\eop

\section{Appendix: about elementary abelian 2-subgroups in the Monster}

\begin{rem}\labtt{about2localsinmonster}
Many examples of 2-local subgroups of the Monster have been known a long time.  A family of maximal 2-locals makes a pretty geometry \cite{rs}.   The next lemma is an elementary argument to prove just what we need in this article.  A classification of maximal 2-locals in the Monster and Baby Monster may be found in preprints at \begin{verbatim}http://www.math.msu.edu/~meier/Preprints/2monster/abstract.html\end{verbatim}
\end{rem}

\begin{lem}\labtt{2Bpure}
Let $G$ be the Monster and let $z$ be a 2B-involution, so that $C:=C_G(z)$ has the form $2^{1+24}Co_1$.    Let $R:=O_2(C)$.   Then, for any $m\in \{1,2,3,4\}$,  $R$ contains a subgroup $E$ so that $E\cap \la z \ra =1$ and
$N_G(E)/C_G(E) \cong GL(m,2)$.
\end{lem}
\pf (Sketch.)
It suffices to prove this for $m=4$.

Let $R \le H \le C$ so that $H/R$ is the centralizer of a 2-central involution in $C/R$, $H/R\cong 2^{1+8}\Omega^+(8,2)$.  Let the coset
$tR$ generate the center of $H/R$.   Then $[R, t]$ is elementary abelian of rank 9 and is a module for $H/R$ such that  the  composition factor $[R,t]/\la z\ra$ has rank 8.  This composition factor may be considered a nonsingular maximal Witt index quadratic space.   In it, take a subspace of rank 4 which is in this way considered totally singular.   Let $S$ satisfy $\la z \ra \le S \le [R,t]$ so that $S/\la z \ra$ is identified with such a subspace.
Let $K:=N_H(S)$, so that the natural map of $K$ to $Aut(S/\la z \ra ) \cong GL(4,2)$ is onto.   Let $E$ be a complement
to $\la z \ra$ in $S$.   The action of $R$ by conjugation on the set of complements
to $\la z \ra$ in $S$ is transitive, whence the factorization
$K=R\, N_H(E)$.   The action of $N_H(E)$ may be identified with the action of $H$ on $S/\la z \ra$.   It follows that $E$ is either $2A$-pure or $2B$-pure.   If $e\in E$ is an involution,
the 2-part of $|C_C(e)|$ is $2^{45}$, whereas the 2-part of the centralizer order for a $2A$-involution is $2^{42}$.  Therefore, $E$ is $2B$-pure.
\eop

\bigskip 

{\bf Acknowledgments.  } This work was started when the first author was a Gehring Visiting Professor in University of Michigan at Ann Arbor in spring 2012. The first author is also supported by NSF grants. The second author acknowledges support from an NSA grant.

\end{document}